\documentclass[11pt]{article}  
\usepackage{amssymb}

\newtheorem{theorem}{Theorem}
\newtheorem{corollary}{Corollary}

\newtheorem{lemma}{Lemma}

\newenvironment{example} {\smallskip\noindent{\bf Example\/}.}{\smallskip\par}
\newenvironment{examples} {\smallskip\noindent{\bf
Examples\/}.}{\smallskip\par}

\newenvironment{remarks} {\smallskip\noindent{\bf Remarks\/}.}{\smallskip\par}
\newenvironment{proof} {\noindent{\em Proof\/}.}{{ $\Box$}\smallskip\par}

\newcommand{\CC}{{\Bbb C}}

\newcommand{\PP}{{\Bbb P}}

\newcommand{\ZZ}{{\Bbb Z}}

\newcommand{\vv}{{\underline v}}
\newcommand{\ttt}{{\underline t}}
\newcommand{\kk}{{\underline k}}
\newcommand{\mm}{{\underline m}}
\newcommand{\1}{{\underline 1}}

\title{On the arc filtration for the singularities of Arnold's lists}
\author{W.Ebeling and S.M.Gusein-Zade
\thanks{Partially supported by the DFG-programme ''Global methods
in complex geometry'' (Eb 102/4--2), grants RFBR--01--01--00739,
INTAS--00-259. Keywords: ring of functions, arcs on a variety,
filtration, Poincar\'e series. AMS Math. Subject Classification:
14B05, 32S10} }

\date{}

\begin{document}

\maketitle

\begin{abstract}
In a previous paper, the authors introduced a filtration on the ring
${\cal O}_{V,0}$ of germs of functions on a germ $(V,0)$ of a complex
analytic variety defined by arcs on the singularity and called the
arc filtration. The Poincar\'e series of this filtration were computed for
simple surface singularities in the 3-space. Here they are computed for
surface singularities from Arnold's lists including uni- and bimodular
ones. The classification of the unimodular singularities by these
Poincar\'e series turns out to be in accordance with their hierarchy
defined by E.~Brieskorn using the adjacency relations. Besides that
we give a general formula for the Poincar\'e series of the arc filtration
for isolated surface singularities which are stabilizations of plane curve
ones.
\end{abstract}

\section*{Introduction}
Let $(V,0)$ be a germ of a complex analytic variety. In \cite{EG}
we defined a filtration on the ring ${\cal O}_{V,0}$ of germs of
functions on the variety $(V,0)$, which we called arc filtration.
An arc $\phi$ on $(V,0)$ is a germ of a complex analytic mapping
$\phi: (\CC,0) \to (V,0)$. For an arc $\phi$ on $(V,0)$ and for
a germ $g \in {\cal O}_{V,0}$ let $v_\phi(g) \in \ZZ_{\geq 0}$ be
the order of the function $g \circ \phi$ at the origin, i.e., the
power of the first non-vanishing term in the power series expansion
$g \circ \phi (\tau) = a \tau ^{v_\phi(g)} + \ldots$\,, $a \neq 0$
(if $g \circ \phi \equiv 0$, then $v_\phi(g)$ is assumed to be equal
to $+\infty$). Let $v(g)$ be the minimum over all arcs $\phi$ on
$(V,0)$ of the orders $v_\phi(g)$. The arc filtration ${\cal O}_{V,0}=F_0
\supset F_1 \supset F_2 \supset \ldots$ on the ring ${\cal O}_{V,0}$
is the filtration by the ideals $F_i= \{ g \in {\cal O}_{V,0} \, :\,
v(g) \geq i\}$. Let $\frak{m}$ be the maximal ideal of the ring
${\cal O}_{V,0}$. Obviously $\frak{m}^i \subset F_i$ and therefore
the ideals $F_i$ have finite codimension in the algebra ${\cal O}_{V,0}$.
Let
$$P_{V,0}(t)= \sum_{i=0}^\infty \dim (F_i/F_{i+1}) \cdot t^i$$
be the Poincar\'e series of the arc filtration on the ring
${\cal O}_{V,0}$.

In \cite{EG}, Poincar\'e series of the arc filtration on the ring
of germs of functions of the rational double points were computed.
A rational double point is an isolated surface singularity defined
by $\{f=0\}$ where $f$ is a simple ($0$-modular) germ of a function
in three variables. For each of these singularities the Poincar\'e
series turns out to coincide with the quasihomogeneous filtration
on the ring of functions of another simple surface singularity in
$(\CC^3,0)$.

In this paper we continue the study of the Poincar\'e series of the
arc filtration of isolated 2-dimensional hypersurface singularities.
We compute the Poincar\'e series $P_{V,0}(t)$ for the following surface
singularities of V.~I.~Arnold's lists \cite{AGV}: the stabilizations
of the curve singularities of multiplicity at most 4 and the
singularities of multiplicity 3 with a reduced 3-jet. These
singularities include all the uni- and bimodular surface singularities.
For the parabolic singularities we have the following table of
Poincar\'e series $P_{V,0}(t)$:
$$
\begin{tabular}{ccc} Notation & Equation & $P_{V,0}(t)$ \\
\hline
$T_{2,3,6}$ ($J_{10}$) & $x^6+y^3+z^2=0$ &
$\frac{1-t^6}{(1-t)(1-t^2)(1-t^3)}$ \\
$T_{2,4,4}$ ($X_9$) & $x^4+y^4+z^2=0$ & $\frac{1-t^4}{(1-t)^2(1-t^2)}$ \\
$T_{3,3,3}$ ($P_8$) & $x^3+y^3+z^3=0$ & $\frac{1-t^3}{(1-t)^3}$ \\
\hline
\end{tabular}
$$
(the notation in brackets is the notation of Arnold \cite{AGV}.)
These Poincar\'e series coincide with the Poincar\'e series of the
natural quasihomogeneous filtrations on the rings of germs of
functions on these singularities. We show that the Poincar\'e
series of the arc filtrations for any uni- or bimodular surface
singularity $(V,0)$ coincides with one of these series. Moreover,
the singularity $(V,0)$ can be deformed to at least one parabolic
singularity \cite{B}. The Poincar\'e series $P_{V,0}(t)$ of the arc
filtration is equal to the corresponding Poincar\'e series of the
parabolic singularity with the minimal Milnor number among those to
which $(V,0)$ can be deformed. This means that the classification
of the unimodular surface singularities by the Poincar\'e series of
their arc filtrations is in accordance with the hierarchy defined by
the adjacency relations described by E.~Brieskorn \cite{B}.

Besides that we give a general formula for the Poincar\'e series
$P_{V,0}(t)$ for an isolated surface singularity $(V,0) \subset
(\CC^3,0)$ which is the stabilization of a (reduced) plane curve
singularity. Moreover, we show that for a hypersurface singularity
with a reduced tangent cone the arc filtration coincides with
the filtration by powers of the maximal ideal.

\section{Stabilizations of curve singularities}
Since a big part of the singularities from Arnold's lists are
stabilizations of plane curve singularities, we start with a
study of this case in general. Let $(C,0)\subset(\CC^2, 0)$ be a germ of a
reduced plane curve given by an equation $f(x,y)=0$ with
$f\in{\cal O}_{\CC^2, 0}$. The stabilization of the plane curve
singularity $(C,0)$ is the surface singularity $(S,0)\subset(\CC^3, 0)$
defined by the equation $f(x,y)-z^2=0$. The surface $(S,0)$ has
an isolated singularity at the origin. Let $\pi:({\cal X, \cal D})
\to(\CC^2, 0)$ be an embedded resolution of the curve singularity
$(C, 0)$. The exceptional divisor ${\cal D}=\pi^{-1}(0)$ of the
resolution $\pi$ is the union $\bigcup_{i=1}^rE_i$ of irreducible
components $E_i$; each of them is a rational curve, i.e., is
isomorphic to the projective line $\CC\PP^1$. For a germ
$g\in{\cal O}_{\CC^2, 0}$, $g\ne0$, let $v_i(g)$ be the multiplicity
of the lifting $\widetilde g=g\circ\pi$ of the function $g$ to the
space $\cal X$ of the resolution $\pi$ along the component $E_i$ of
the exceptional divisor $\cal D$. We shall call $m_i=v_i(f)$ the
multiplicity of the component $E_i$ of the exceptional divisor ($f$
is the left hand side of the equation of the curve $(C,0)$). Let
$\sigma_i=1$ (respectively $\sigma_i=2$) if the multiplicity $m_i$
is even (respectively odd). Suppose that the resolution $\pi$ has the
property that the strict transform $\widetilde C$ of the curve $C$
intersects the exceptional divisor $D$ only at points of the components
$E_i$ with even multiplicity $m_i$. Such a resolution can be obtained
from an arbitrary (say, the minimal) one by blowing up intersection
points of the strict transform $\widetilde C$ with the components $E_i$
with odd multiplicity $m_i$.

The resolution $\pi$ defines a multi-index filtration on the ring
${\cal O}_{\CC^2, 0}$ of germs of functions of two variables which
is called {\em divisorial} and which is defined by the valuations
$v_i(\cdot)$, i.e., the corresponding ideal $J(\vv)$
($\vv=(v_1, \ldots, v_r) \in \ZZ^r$) is equal to
$\{g\in {\cal O}_{\CC^2, 0}: v_i(g)\ge v_i,\, i=1, \ldots, r \}$. (Pay
attention that we defined the ideal $J(\vv)$ for all $\vv \in \ZZ^r$,
not only for those with non-negative components). Let
$P_\pi(t_1, \ldots, t_r)$ be the Poincar\'e series of this filtration
(see, e.g., \cite{DG}). The Poincar\'e series $P_\pi(t_1, \ldots, t_r)$
is a power series in $t_1$, \dots, $t_r$ defined as
$$
\frac{\left(
\sum_{\vv \in \ZZ^r}\dim\left(J(\vv)/J(\vv+\1)\right)\ttt^\vv
\right)
\cdot\prod_{i=1}^r(t_i-1)}{t_1\cdot\ldots\cdot t_r-1}.
$$
Here $\ttt^\vv= t_1^{v_1}\cdot\ldots\cdot t_r^{v_r}$, $\1=(1, \ldots,
1)\in\ZZ^r$,
$\sum_{\vv \in \ZZ^r}\dim\left(J(\vv)/J(\vv+\1)\right)\ttt^\vv$ is
a formal Laurent series in the variables $t_1$, \dots, $t_r$ infinite
in all directions for $r\ge 2$ (for $r=1$ it is a power series).

For a component $E_i$ of the exceptional divisor $\cal D$ of the
resolution $\pi$, let ${\stackrel{\circ}{E}}_i$ be its "smooth part",
i.e., $E_i$ without the intersection points with other components of
the exceptional divisor $\cal D$, let $\widetilde L_i$ be a germ of a
non-singular curve transverse to the component $E_i$ at a smooth point,
i.e., at a point of ${\stackrel{\circ}{E}}_i$. Let the
blow-down $L_i=\pi(\widetilde L_i)$ of the curve $\widetilde L_i$
be given by an equation $g_i=0$, $g_i\in {\cal O}_{\CC^2, 0}$, and let
$m_{ij}:=v_j(g_i)$. The
matrix $(m_{ij})$ is symmetric and is the inverse to minus the matrix of
intersections of the components $E_i$.
Let $\mm_i:=(m_{i1}, \ldots, m_{ir})$. In \cite{DG} it was shown that
$$
P_\pi(t_1, \ldots, t_r)=
\prod_{i=1}^r(1-\ttt^{\mm_i})^{-\chi({\stackrel{\circ}{E}}_i)},
$$
where $\chi(A)$ is the Euler characteristic of a space $A$.

For a power series $P(t_1, \ldots, t_r)$, the {\em reduction}
$\overline{P}(t)$ of this series is
the series in one variable $t$ obtained from $P$ by substituting each
monomial $\ttt^\kk$,
$\kk=(k_1, \ldots, k_r)\in\ZZ_{\ge0}^r$, by the monomial $t^k$ with
$k=\min\limits_{1\le i\le r}
k_i$ \cite{EG}. For example, for $P(t_1, t_2)=1+t_1t_2^2+t_1^2t_2+t_1^2t_2^2$,
one has
$\overline{P}(t)=1+2t+t^2$.

\begin{theorem}\label{Thm1} The Poincar\'e series $P_{S,0}(t)$ of the arc
filtration of the
stabilization $(S,0)$ of the curve $(C,0)$ is equal to the reduction
$\overline{Q}(t)$ of the
series
$$ Q(t_1, \ldots, t_r)=(1+t_1^{\sigma_1 m_1/2}\cdot\ldots\cdot
t_r^{\sigma_r m_r/2})\cdot
P_\pi(t_1^{\sigma_1}, \ldots, t_r^{\sigma_r})\, .
$$
\end{theorem}

\begin{remarks} {\bf 1.} In $\cite{EG}$, there is given a somewhat
similar formula for the Poincar\'e series $P_{V,0}(t)$ for an
arbitrary irreducible isolated singularity $(V,0)$ in terms of its
resolution. However, for an arbitrary singularity (in contrast to a
plane curve one) sometimes it is difficult to construct a resolution
and, moreover, no general formula for the Poincar\'e series of the
corresponding set of divisorial valuations is known.\newline
{\bf 2.} Though the series $Q(t_1, \ldots, t_r)$ depends on $r$
variables, $r$ being the number of irreducible components of the
exceptional divisor $\cal D$, the minimal exponents usually can be met
only at variables corresponding to a few first components of the
exceptional divisor (in the sense of the order of there appearance
in the process of resolution by blow-ups). Often this reduces
the calculations considerably.
\end{remarks}

\begin{example} The rational double point $D_4$ is the stabilization of the
curve singularity
$x^3+y^3=0$. This curve is resolved after blowing up once. However, the
exceptional divisor has
odd multiplicity (equal to 3). Therefore, to get a suitable resolution, one
has to blow up the 3
intersection points of the strict transform of the curve with the
exceptional divisor. This leads
to the resolution with the exceptional divisor consisting of 4 components
with the intersection
matrix
$$
\left( \begin{array}{cccc} -4 & 1 & 1 & 1 \\ 1 & -1 & 0 & 0 \\ 1 & 0 & -1 &
0 \\1 & 0 & 0 & -1
\end{array} \right).
$$
Therefore
$$
(m_{ij}) = \left( \begin{array}{cccc} 1 & 1 & 1 & 1 \\1 & 2 & 1 & 1\\1 & 1
& 2 & 1\\1 & 1 & 1 & 2
\end{array} \right)
$$
and hence
$$
P_\pi(t_1,t_2,t_3,t_4)=
\frac{1-t_1t_2t_3t_4}{(1-t_1t_2^2t_3t_4)(1-t_1t_2t_3^2t_4)(1-t_1t_2t_3t_4^2)}.
$$
Theorem~\ref{Thm1} implies that the Poincar\'e series of the arc filtration is
the reduction
$\overline{Q}(t)$ of the series
$$
(1+t_1^3t_2^2t_3^2t_4^2)
\frac{1-t_1^2t_2t_3t_4}{(1-t_1^2t_2^2t_3t_4)(1-t_1^2t_2t_3^2t_4)(1-t_1^2t_2t_3t_
4^2)}
$$
which is equal to
$$
\frac{1-t^3}{(1-t)^2(1-t^2)}.
$$
In \cite{EG} this result was obtained using a general description
of the Poincar\'e series in several variables for rational surface
singularities from \cite{CDG}.
\end{example}

The ring ${\cal O}_{S, 0}$ of germs of functions on the stabilization
$(S,0)$ is a free ${\cal
O}_{\CC^2, 0}$--module of rank 2. Moreover each germ $H\in{\cal O}_{S, 0}$
can in a
unique way be written in the form $H=h_1+z\cdot h_2$ with $h_i\in {\cal
O}_{\CC^2, 0}$,
$i=1,2$. Thus one can identify ${\cal O}_{S, 0}$ with ${\cal O}_{\CC^2,
0}\oplus
z\cdot{\cal O}_{\CC^2, 0}$. The arc filtration on the ring ${\cal O}_{S,
0}$ defines
corresponding filtrations on the summands ${\cal O}_{\CC^2, 0}$ and
$z\cdot{\cal
O}_{\CC^2, 0}$.

\begin{lemma} \label{Lemma1}
One has $v(H)=\min\left(v(h_1), v(z\cdot h_2)\right)$. Therefore the
Poincar\'e series
$P_{S,0}(t)$ is the sum of the respective Poincar\'e series of these
filtrations on
${\cal O}_{\CC^2, 0}$ and $z\cdot{\cal O}_{\CC^2, 0}$.
\end{lemma}

\begin{proof} Obviously $v(H)\ge \min\left(v(h_1), v(z\cdot h_2)\right)$.
Suppose that $v(H)>
v:=\min\left(v(h_1), v(z\cdot h_2)\right)$. This implies that there exists
an arc $\psi:(\CC, 0)\to(S,
0)$ on $S$ such that $h_1\circ\psi(\tau)=a\cdot\tau^v+\ldots\,$,
$(z\cdot h_2)\circ\psi(\tau)=-a\cdot\tau^v+\ldots\,$, $a\ne0$. Consider the
arc $\psi'$ defined by
$\psi'(\tau)=\iota\circ\psi(\tau)$ where $\iota$ is the natural involution
on the surface $S$:
$(x, y, z)\mapsto(x, y, -z)$. Then $v_{\psi'}(H)=v$, a contradiction. The
statement about the
Poincar\'e series follows from the fact that $\{H\in{\cal O}_{S,
0}:\,v(H)\ge i\}=\{h_1\in{\cal
O}_{\CC^2, 0}:\,v(h_1)\ge i\} \oplus \{h_2\in{\cal O}_{\CC^2, 0}:\,v(z\cdot
h_2)\ge i\}$.
\end{proof}

For $h\in {\cal O}_{\CC^2, 0}$, let $\widetilde v_i(h):= \sigma_i\cdot
v_i(h)$. The functions
$\widetilde v_i(\cdot)$ are not valuations (they do not respect the
multiplication), however,
they define a multi-index filtration on the space ${\cal O}_{\CC^2, 0}$ as
well:
$\widetilde J(\vv)=\{h\in {\cal O}_{\CC^2, 0}: \widetilde v_i(h)\ge v_i,\,
i=1, \ldots, r \}$ for
$\vv\in\ZZ^r$.

Let $\widetilde{v}_i(z\cdot h):= \sigma_i\cdot(v_i(h)+(m_i/2))$. Here, with
some abuse of notation, one can
say that $v_i(h)+(m_i/2)$ is the order along the component $E_i$ of the
lifting to the space of the
resolution of the function $\sqrt{f}\cdot h$. Then the functions
$\widetilde{v}_i$ define a multi-index
filtration on the space $z\cdot{\cal O}_{\CC^2,0}$ as well.

\begin{lemma} \label{Lemma2} One has
$$
v(h)=\min_{1 \leq i \leq r} \widetilde{v}_i(h),
\qquad
v(z\cdot h)=\min_{1 \leq i \leq r} \widetilde{v}_i(z\cdot h).
$$
\end{lemma}

\begin{proof}
Under the projection $(x,y,z) \mapsto (x,y)$, an arc on the surface
$(S,0)$ is mapped to a plane arc. Moreover, each arc in $(\CC^2,0)$
is the projection of one or two geometrically different arcs on the
surface $(S,0)$. The number of arcs on $(S,0)$ which are mapped to
the same plane arc $\phi: (\CC,0) \to (\CC^2,0)$ depends on the order
$v_\phi(f)$ of the function $f$ on the arc $\phi$ and is equal to 1
or 2 if this order is odd or even (and finite) respectively. Moreover,
in the first case the arc on $(S, 0)$ is a two-fold covering of its
projection, in the second case the projection defines an isomorphism
between the curve in the plane and each of its preimages. If the
order $v_\phi(f)$ is infinite, i.e., if the arc $\phi$ lies in the
curve $C=\{f=0\}$, the arc $\phi$ has one preimage isomorphic to it.

Let us prove the first statement of the lemma. Let
$\widetilde{S} \subset {\cal X} \times \CC$ be defined by the equation
$\widetilde{f}-z^2=0$ where $\widetilde{f}=f \circ \pi$ is the lifting
of the function $f$ to the space $\cal X$ of the resolution. The surface
$\widetilde{S}$ is a modification of the surface $S$. Arcs on $S$
arriving at the origin are in one-to-one correspondence with arcs on
$\widetilde{S}$ arriving at the exceptional divisor ${\cal D} \times
\{0\}$. If we take an arc on $\cal X$ which intersects the exceptional
divisor ${\cal D}$ at a generic point of the component $E_i$ then the
order of the function $\widetilde{h}$ along a preimage of this arc in
$\widetilde{S}$ is just equal to $\widetilde{v}_i(h)$. (Note that here
we use the fact that the strict transform of the curve $\{f=0\}$
intersects only the components $E_i$ of the exceptional divisor ${\cal D}$
with even multiplicity $m_i$.) Therefore $v(h) \leq \min_{1 \leq i \leq r}
\widetilde{v}_i(h)$.

Let us prove the opposite inequality. Suppose that an arc $\phi$
on $\cal X$ intersects the exceptional divisor ${\cal D}$ at a point
of $\stackrel{\circ}{E}_i$ with multiplicity $\mu_i \geq 2$. If both
$\mu_i$ and $m_i$ are odd then the order of the function $\widetilde{h}$
along a preimage of the arc $\phi$ in $\widetilde{S}$ is greater than
or equal to $2\mu_iv_i(h)$ (since the preimage is a two-fold covering
of the arc $\phi$). If at least one of the numbers $\mu_i$ and $m_i$
is even then this order is greater than or equal to $\mu_i v_i(h)$.
In any case it is not less than $\widetilde{v}_i(h)$.

Suppose that an arc $\phi$ on $\cal X$ intersects the exceptional
divisor ${\cal D}$ at the intersection point of two divisors $E_i$
and $E_j$. After several additional blow-ups of the intersection
points of the components of the exceptional divisor, the strict transform
of the arc $\phi$ will intersect a new component at a smooth point.
Therefore it is sufficient to prove that an additional blow-up of
the intersection points of two components (say, $E_i$ and $E_j$) does not
change the number $\min_{\{i\}} \widetilde{v}_i(h)$. Let the new component
created by blowing up be $E_\alpha$. Then $v_\alpha(h)\ge v_i(h)+v_j(h)$,
$m_\alpha=m_i+m_j$. It is not difficult to show that
$\widetilde{v}_\alpha(h)\ge\min\{\widetilde{v}_i(h),\widetilde{v}_j(h)\}$
(this can be verified, e.g., considering three cases: 1) both $m_i$ and
$m_j$ are even; 2) both $m_i$ and $m_j$ are odd; 3) $m_i$ is even,
$m_j$ is odd). This implies the statement.

The proof of the other statement is analogous.
\end{proof}

\begin{lemma} \label{Lemma3}
The Poincar\'e series of the multi-index filtration on ${\cal O}_{\CC^2,0}$
(respectively on $z\cdot{\cal O}_{\CC^2,0}$) defined by the functions
$\widetilde{v}_i( \cdot)$ is equal to
$P_\pi(t_1^{\sigma_1}, \ldots , t_r^{\sigma_r})$ (respectively to
$t_1^{\sigma_1m_1/2} \cdots t_r^{\sigma_r m_r/2}
P_\pi(t_1^{\sigma_1}, \ldots , t_r^{\sigma_r})$.)
\end{lemma}

\begin{proof}
A convenient way to see this is to write these Poincar\'e series
as integrals with respect to the Euler characteristic (see, e.g.,
\cite{DG}). For example, the Poincar\'e series of the multi-index
filtration defined by the functions $\widetilde{v}_i( \cdot )$ on
the space $z\cdot{\cal O}_{\CC^2,0}$ is equal to
$$
\int_{\PP{\cal O}_{\CC^2,0}}
t_1^{\widetilde{v}_1(zh)} \cdots t_r^{\widetilde{v}_r(zh)} d\chi.
$$
Comparison of these integrals gives the statement.
\end{proof}

\noindent {\em Proof of Theorem~\ref{Thm1}.}
The Poincar\'e series of the one-index filtration defined by the
function $\min_{1 \leq i \leq r} \widetilde{v}_i( \cdot)$ is the
reduction of the (multi-variable) Poincar\'e series of the multi-index
filtration defined by the functions $\widetilde{v}_i(\cdot)$. (This
can be derived using the same arguments as in \cite[Proposition~2]{EG}.)
Now the statement follows from Lemmas 1--3.
$\Box$

\addvspace{3mm}

As a corollary of Theorem~\ref{Thm1} one has the following statement.
For a germ $f \in {\cal O}_{\CC^2,0}$ let $m(f)$ be its multiplicity,
i.e., $f \in \frak{m}^{m(f)}\setminus \frak{m}^{m(f)+1}$. Suppose that
an irreducible curve is tangent to a smooth one. In local coordinates
$(u,v)$ with the $u$-axis representing the smooth curve the first one
can be given by an equation $v = p(u)$ where $p(u)$ is a Puiseux series.
We call the first non-trivial exponent of this series the {\em order of
contact} between the curves.

\begin{corollary} \label{Cor1} {\rm (i)} If $m(f)$ is even then
$$P_{S,0}(t)=\frac{1-t^m}{(1-t^{m/2})(1-t)^2}.$$

\noindent{\rm (ii)} Let $m(f)$ be odd and let all components of
the curve $C$ be tangent to the same smooth curve with order of
contact at least 2. Then
$$
P_{S,0}(t)=\frac{1-t^{2m}}{(1-t^m)(1-t^2)(1-t)}.
$$
\end{corollary}

\begin{proof}
Let us number the components of the exceptional divisor in the
order of their appearance in the process of the resolution by
blow-ups. If $m(f)$ is even the statement follows from the
fact that $m_1=m(f)$, $\sigma_1=1$, and for any $i,j$ one has
$m_{i1}=1 \leq m_{ij}$. If $m(f)$ is odd then $m_1=m(f)$, $m_2=2m(f)$,
$m_{i1}=1$, $m_{i2}=2$ for $i\geq 2$ and $m_{ij} \geq m_{i2}$ for
$j \geq 2$. This implies the statement.
\end{proof}

\section{Singularities with a reduced tangent cone}
Let a hypersurface singularity $(V,0) \subset (\CC^{n+1},0)$ be
given by an equation $f=0$ with $f=f_d + f_{d+1} + \ldots$ where
$f_d$ is the homogeneous part of $f$ of degree $d$, $f_d\ne 0$.
The tangent cone to the singularity $(V,0)$ is given by the
equation $f_d=0$.

\begin{theorem}\label{Thm2}
If the tangent cone is reduced then the arc filtration on the ring
${\cal O}_{V,0}$ coincides with the filtration by powers of the
maximal ideal. Therefore
$$
P_{V,0}(t)= \frac{1-t^d}{(1-t)^{n+1}}.
$$
\end{theorem}

\begin{proof}
Obviously $\frak{m}^i \subset F_i$. To prove that $\frak{m}^i = F_i$
we shall show that for $h \in \frak{m}^i \setminus \frak{m}^{i+1}$
we have $v(h)=i$. This follows from the fact that there are many
smooth arcs on the hypersurface $V$. The fact that
$h \in \frak{m}^i \setminus \frak{m}^{i+1}$ means that a representative
of $h$ in ${\cal O}_{\CC^{n+1},0}$ (which we also denote by $h$)
can be chosen of the form $h=h_i + h_{i+1} + \ldots$ where the
homogeneous part $h_i$ of degree $i$ is not divisible by $f_d$.
Therefore there is a regular point $p$ in the hypersurface
$\{f_d=0\} \subset \CC\PP^n$ such that $h_i(p) \neq 0$. After
blowing up the origin in $\CC^{n+1}$, the strict transform
$\widetilde{V}$ of the hypersurface $V$ is non-singular at the
point $p$ of the exceptional divisor $\CC\PP^n$ and intersects
$\CC\PP^n$ at this point transversely. An arc in the strict
transform $\widetilde{V}$ transverse to $\CC\PP^n$ blows down
to a smooth arc along which the order of the function $h$ is
equal to $i$. Therefore $v(h) \leq i$.
\end{proof}

\begin{examples}
{\bf 1.} In particular, Theorem~\ref{Thm2} implies the following statement:
if a hypersurface singularity $(V,0)$ is the double suspension of a
hypersurface singularity $\{f=0\} \subset \CC^n$ with $f \in \frak{m}^2$
(i.e., it is the hypersurface singularity in $\CC^{n+2}$ defined by
the equation $f(z_1, \ldots , z_n)+x^2+y^2=0$) then
$$
P_{V,0}(t)=\frac{1-t^2}{(1-t)^{n+2}}\quad \mbox{\cite[Proposition~2]{EG}}.
$$

\noindent{\bf 2.} If $f(x,y,z)=xyz + terms\ of\ higher\ degree$ then,
for the
surface singularity $V=\{f=0\}$, one has
$$
P_{V,0}(t)=\frac{1-t^3}{(1-t)^3}.
$$
\end{examples}

\section{Singularities of Arnold's lists}
\begin{theorem}
The Poincar\'e series $P_{V,0}(t)$ of the arc filtration are given
by the following list:

{\rm (i)} For the stabilizations of the curve singularities with non-zero
3-jet except the simple
ones, i.e.,
$J_{k,i}$ ($i \geq 0$),
$E_{6k}$, $E_{6k+1}$, and $E_{6k+2}$ (in all cases $k \geq 2$), we have
$$P_{V,0}(t)=\frac{1-t^6}{(1-t)(1-t^2)(1-t^3)}.$$

{\rm (ii)} For the stabilizations of the curve singularities with zero
3-jet and non-zero 4-jet,
i.e., for those of classes ${\bf X}$, ${\bf Y}$,
${\bf Z}$, and ${\bf W}$, we have
$$P_{V,0}(t)=\frac{1-t^4}{(1-t)^2(1-t^2)}.$$

{\rm (iii)} For the singularities of multiplicity 3 with a reduced 3-jet,
i.e., for those of the
series ${\bf Q}$, ${\bf S}$, and ${\bf U}$, we have
$$P_{V,0}(t)=\frac{1-t^3}{(1-t)^3}.$$
\end{theorem}

\begin{proof} Statement (i) is a particular case of Corollary~\ref{Cor1}
with $m=3$. Statement
(ii) is a particular case of this corollary with $m=4$. Finally, (iii) is a
particular case of
Theorem~\ref{Thm2}.
\end{proof}

\bigskip
\noindent Universit\"{a}t Hannover, Institut f\"{u}r Mathematik \\
Postfach 6009, D-30060 Hannover, Germany \\
E-mail: ebeling@math.uni-hannover.de\\

\medskip
\noindent Moscow State University, Faculty of Mechanics and Mathematics\\
Moscow, 119992, Russia\\
E-mail: sabir@mccme.ru

\end{document}